\documentclass[12pt,twoside,a4paper]{amsart}
\usepackage{amsmath,amssymb,amscd,amsfonts,latexsym,a4wide}
\newtheorem{thm}{Theorem}[section]

\newtheorem{prop}[thm]{Proposition}

\theoremstyle{definition}

\newtheorem{prg}[thm]{}

\theoremstyle{remark}
\newtheorem{remark}[thm]{Remark}

\numberwithin{equation}{section}

\newcommand{\thmref}[1]{Theorem~\ref{#1}}

\newcommand{\propref}[1]{Proposition~\ref{#1}}

\newcommand{\formref}[1]{(\ref{#1})}
\newcommand{\prgref}[1]{{\bf\ref{#1}}}

\newcommand{\Cc}{{\mathbb C}}
\newcommand{\Aa}{{\mathbb A}}
\newcommand{\Zz}{{\mathbb Z}}
\newcommand{\Qq}{{\mathbb Q}}
\newcommand{\Pp}{{\mathbb P}}
\newcommand{\Rr}{{\mathbb R}}

\newcommand{\Ff}{{\mathbb F}}

\newcommand{\Ob}{{\mathcal O}}
\newcommand{\Lb}{{\mathcal L}}
\newcommand{\Hom}{\operatorname{Hom}}

\renewcommand{\leq}{\leqslant}
\renewcommand{\geq}{\geqslant}

\def\rep#1{\bysame}
\def\bib[#1]#2<#3>#4|#5(#6)#7-#8.{\bibitem{#1}
{\sc #2},\ #3,\ {\it #4},\ {\bf #5}\ (#6)\ #7--#8.}
\def\tbib[#1]#2<#3>#4.{\bibitem{#1} {\sc #2},\ {\it #3},\ #4.}
\def\bibt[#1]#2<#3>#4|#5(#6)#7-#8{\bibitem{#1}
{\sc #2},\ #3,\ in\ {\it #4}\ ed.\ {#5},\ (#6)\ #7--#8.}
\def\bibit[#1]#2<#3>{\bibitem{#1} {\sc#2},\ #3}
\def\prebib[#1]#2<#3>{\bibitem{#1} {\sc #2},\ {\it #3}, Preprint.}
\def\toabib[#1]#2<#3>{\bibitem{#1} {\sc #2},\ #3,\ to appear.}

\title[Mumford's Degree of Contact and Diophantine
Approximations]{Mumford's Degree of Contact and
Diophantine Approximations}
\author[Roberto G. Ferretti]{Roberto G. Ferretti}
\address{Insitute for Advanced Study, School of Mathematics,
Olden Lane, Princton, NJ 08540}
\address{Institut des Hautes \'Etudes Scientifiques, 35 Route de
Chartres, 91440 Bures-sur-Yvette, France}
\email{ferretti@ias.edu, ferretti@ihes.fr}
\date{\today}
\begin{document}

\maketitle
\section*{Introduction} The Schmidt Subspace Theorem
establishes that the solutions of some particular systems
of diophantine approximations in projective spaces
accumulates on a finite number of linear subspaces
(\cite{Sc}, Theorem
$V.1.D'$). One may state the following question:
given a subvariety $X$ of a projective space $\Pp^n$,
does there exists a system of diophantine approximations
on $\Pp^n$ whose solutions are Zariski dense in $\Pp$ but
lie in finitely many proper subvarieties of $X$? One can gain insight
into this problem using a theorem
of G.~Faltings and G.~W\"ustholz (\cite{FW}, Theorem $7.3$, or
\cite{Fe} for a quantitative version). Their construction requires
the hypothesis that the sum of some expected
values (see \formref{uno.due}) has to be large. This
sum turns out to be proportional to a degree of contact of a
weighted flag of sections over the variety $X$ (see
\formref{uno.cinque}). The degree of contact was first introduced
by D.~Mumford (see for example \cite{Mu}) as a tool for
the construction of moduli spaces for families of algebraic
schemes (\cite{Mu} for
curves, \cite{Gi} for canonically embedded minimal
surfaces, or \cite{Mo} for ruled surfaces). Indeed, it measures the
semistability of the Chow (or Hilbert) point of $X$ under the action
of an appropriate reductive algebraic group. Whence, the lower bound
in the Faltings-W\"ustholz theorem may be translated into a geometric
invariant theoretic language. Moreover, the question of understanding
when a system of diophantine approximations on $X$ is not under control
of the Schmidt Subspace Theorem relies on the Chow-unstbility of $X^s$, for
some sufficiently large $s>0$. Unstability
is usually associated to singular varieties. We will see
this for some unstable local rings and for a family of
elliptic surfaces. However ruled surfaces and some
blow-ups give examples of non singular varieties having
an unstable embedding. In general it is quite difficult
to find examples of unstable projective varieties. All noteworthy instances
we found have a common style: the destabilizing flags are constructed
by the vanishing of sections along special divisors. Moreover, in
order to control all the dimensions involved, it is crucial to use
Riemann-Roch, as well as some result on the vanishing of higher cohomologies.\\
{\small It is a pleasure to thank J.-B.~Bost, M.~McQuillan, C.~Soul\'e, G.~W\"ustholz for several discussions and remarks; IAS and IHES for the ospitality;
and Swiss National Science Foundation and NSF with the grant DMS 9304580 for
support}.
\section{Preliminaries}
\begin{prg} Let $K$ be an algebraic number field. Denote its
ring of integers by ${\mathcal O}_K$ and its collection
of places (equivalence classes of absolute values) by
$M_K$. For $v\in M_K$, $x\in K$, we define the absolute
value $|x|_v$ by
\begin{enumerate}
\item $|x|_v=|\sigma(x)|^{1/[K:\Qq]}$ if $v$ corresponds to the
embedding $\sigma:K\hookrightarrow\Rr$;
\item $|x|_v=|\sigma(x)|^{2/[K:\Qq]}=|\overline{\sigma}(x)
|^{2/[K:\Qq]}$ if $v$ corresponds to the the pair of
conjugate embeddings
$\sigma,\overline{\sigma}:K\hookrightarrow\Cc$;
\item $|x|_v=(N{\mathfrak p})^{-\text{ord}_{\mathfrak p}(x)/[K:\Qq]}$
if $v$ corresponds to the prime ideal ${\mathfrak p}$ of
${\mathcal O}_K$.
\end{enumerate}
Here $N{\mathfrak p}=\#({\mathcal O}_K\not {\mathfrak p}$
is the norm of ${\mathfrak p}$ and $\text{ord}_{\mathfrak
p}(x)$ the exponent of ${\mathfrak p}$ in the prime ideal
decomposition of $(x)$, with $\text{ord}_{\mathfrak
p}(0):=\infty$. In case $(i)$ or $(ii)$ we call $v$ real
infinite or complex infinite, respectively and write
$v|\infty$; in case $(iii)$ we call $v$ finite and write
$v\nmid\infty$. These absolute values satisfy the {\it
Product Formula}
\begin{eqnarray}
\label{uno.sei}
\prod_{v\in M_K}|x|_v=1,\ \text{for }x\in K^*.
\end{eqnarray}
The height of ${\bf x}=(x_0,\cdots,x_N)\in K^{N+1}$ with
${\bf x}\neq 0$ is defined as follows: for $v\in M_K$ put
\begin{eqnarray*}
|{\bf
x}|_v=(\sum_{i=0}^N|x_i|_v^{2[K:\Qq]})^{1/2[K:\Qq]}
\text{if }v\text{ is real infinite},\\
|{\bf x}|_v=(\sum_{i=0}^N|x_i|_v^{[K:\Qq]})^{1/[K:\Qq]}\
\text{if }v\text{ is complex infinite},\\
|{\bf x}|_v=\max\{|x_0|_v,\cdots,|x_N|_v\}\ \text{if
}v\text{ is finite}.
\end{eqnarray*}
Now define
$$
h({\bf x})=h(x_0,\cdots,x_N)=\sum_{v\in M_K}\log|{\bf
x}|_v.
$$
By the Product Formula \formref{uno.sei} this defines a
function on $\Pp^N(K)$. Further, $h({\bf x})$ depends
only on ${\bf x}$ and not on the choice of the number
field $K$ containing the coordinates of ${\bf x}$, in
other words the function $h({\bf x})$ extends to a
function on $\Pp^N(\overline{\Qq})$, where
$\overline{\Qq}$ is the algebraic closure of $\Qq$ inside
$\Cc$.\\ It is possible to define the height $h(X)$ of a
projective variety $X\subseteq
\Pp^N$ defned over $K$ (\cite{SABK} for a good reference). This invariant
is comparable to the height of the coefficients of the
(Cayley-Bertini-van der Waerden-)Chow point of $X$.
\end{prg}

\begin{prg} Let $E$ be a vector space of rank $N+1$ over the
number field $K$ and $E^\vee=\Hom (E,K)$. Consider a
closed subvariety $X\subseteq \Pp(E^\vee)$ in the
projective space of lines of $E^\vee$, of dimension $d$.
We choose a basis $l_0,\cdots ,l_N\in E$,  fix non
negative real numbers $r_0\geq \cdots\geq r_N$ and let
${\bf r}=(r_0,\cdots ,r_N)$. When $m$ is large enough,
$m\geq m_0$ say, the cup product map
$$
\varphi:E^{\otimes m}\to H^0(X,{\mathcal O}(m))
$$
is surjective, so that $H^0(X, {\mathcal O}(m))$ is
generated be the monomials
$$
l_0^{\alpha_0}\cdots
l_N^{\alpha_N}=\varphi(l_0^{\otimes\alpha_0}\otimes\cdots
\otimes l_N^{\otimes\alpha_N}),
$$
with $\alpha_0+\cdots+\alpha_N=m$. A special basis is a
basis of $H^0(X,{\mathcal O}(m))$ made of such elements.
We define the weight of $l_i$ to be $r_i$,
$i=0,\cdots,N$, the weight of a monomial in $E^{\otimes
m}$ to be the sum of the weights of the $l_i$'s occurring
in it, and the weight of a monomial $u\in H^0(X,{\mathcal
O}(m))$ to be the minimum $w_{\bf r}(u)$ of the weights
of the monomials in the $l_i$'s mapping to $u$ by
$\varphi$. The weight of a special basis is the sum of of
the weights of its elements. Finally, $w_{\bf r}(m)$
denotes the minimal weight among all special bases of
$H^0(X,{\mathcal O}(m))$. Fix a special basis with
minimal weight and denote by $w_1,\cdots w_M$ the weights
of its elements in increasing order. Let $F^j$ be the
subspace of $H^0(X,{\mathcal O}(m))$, $j=1,\cdots,M$,
generated by those monomials $u$ for which the weight
$w_{\bf r}(u)$ is smaller than $w_j$. We define a
probability measure on $\Rr$ through the density
function:
\begin{eqnarray}
\label{uno.uno}
\rho_m(x)=\frac{1}{h^0(X,{\mathcal O}(m))}
\sum_{w_j\leq x\cdot m}\dim (F^{j+1}/F^j)=\min_{w_j\leq x\cdot m}
\frac{\dim F^j}{h^0(X,{\mathcal O}(m))}.
\end{eqnarray}
Then the expected value $E(\rho_m)$ of the probability
measure is
\begin{eqnarray}
\label{uno.due}
E(\rho_m)=\frac{w_{\bf r}(m)}{m\cdot h^0(X,{\mathcal
O}(m))}.
\end{eqnarray}
\end{prg}

\begin{prg} \label{prg:quattro} For integers
$r_0\geq \cdots\geq r_N$ and a basis $l_0,\cdots,l_N$ of
$E$, the {\it degree of contact} to $X$ is the integer
$e_{\bf r}(X)$, such that when $m$ goes to infinity
\begin{eqnarray}
\label{uno.tre}
w_{\bf r}(m)=e_{\bf r}(X)\frac{m^{d+1}}{(d+1)!}+O(m^d),
\end{eqnarray}
(\cite{Mo}, Corollary $3.3$). By the usual theory of
Hilbert polynomials
\begin{eqnarray}
\label{uno.quattro}
h^0(X,{\mathcal O}(m))=\deg(X)\frac{m^d}{d!}+O(m^{d-1}).
\end{eqnarray}
Combining \formref{uno.tre} with \formref{uno.quattro} we
obtain
\begin{eqnarray}
\label{uno.cinque}
E_\infty(X):=\lim_{m\to\infty}E(\rho_m)=\frac{e_{\bf
r}(X)}{(d+1)\deg(X)}.
\end{eqnarray}
\end{prg}

\begin{prg} Let $L\supseteq K$ a finite field extension
and $\Sigma$ be a finite set of places of $L$ containig
all infinite places. For each $v\in \Sigma$ we choose a
basis $l_{v,0},\cdots,l_{v,N}$ of $E\otimes_KL$ and non
negative real numbers $r_{v,0}\geq \cdots\geq r_{v,N}$ as
above. This defines, for all big enough integers $m$, a
probability measure
\formref{uno.uno} and consequently a real number
$E_{v,\infty}(X\otimes_KL)$ as in \formref{uno.cinque}.
\begin{thm}
\label{th:uno}
Let us suppose that for some real number $0<\delta<1$ we
have,
\begin{eqnarray}
\label{stimaexl}
\sum_{s\in \Sigma}E_{v,\infty}(X\otimes_KL)=1+\delta.
\end{eqnarray}
Then there are two effectively computable constants $c_1,
c_2$ depending only on $K$, $\Sigma$, $\delta$, $N$,
$\deg(X)$, $h(X)$, and the set $\{l_{v,i};\ v\in \Sigma,\
i=0,\cdots,N\}$ such that all points ${\bf x}\in X(K)$
with
\begin{eqnarray}
\label{uno.sette}
h({\bf x})>c_1,
\end{eqnarray}
and
\begin{eqnarray}
\label{uno.otto}
\log(\frac{|l_{v,i}({\bf x})|_v}{|{\bf x}|_v})\leq -r_{v,i}h({\bf x})\quad
v\in \Sigma,\ i=0,\cdots ,N,
\end{eqnarray}
are contained in at most $c_2$ proper subvarieties of $X$
of degree not exceeding $\deg(X)$.
\end{thm}
\proof \cite{Fe} Theorem $10.2$  and \cite{FW} Theorem $7.3$.\qed
\begin{remark}
\label{re:boh} For this theorem it is allowed that all weights $r_{v,i}$
are rationals, but \formref{stimaexl} still holds. This is no restriction
since everything can be multiplied by a common denominator. We will tacitly
assume this in the next sections.
\end{remark}
\end{prg}

\begin{prg}\label{pr:tre} There is an intersection theoretic formula 
expressing
the degree of contact as the degree of divisor on a
suitable modification of $X$. \\ Let ${\mathcal
O}_{X\times\Aa^1}(1)={\mathcal O}_X(1)\otimes{\mathcal
O}_{\Aa^1}$ and $t$ the coordinate of $\Aa^1$. For integers
$r_0\geq \cdots\geq r_N$ and a basis $l_0,\cdots,l_N$ of $E$
we associate the $K[t]$-submodule $I$
of $H^0(X\times\Aa^1,{\mathcal O}_{X\times\Aa^1}(1))$
generated by $\{t^{r_i}l_i,\ i=0,\cdots,N\}$ and an ideal
sheaf $J\subseteq {\mathcal O}_{X\times\Aa^1}$ defined by
$$
J\cdot{\mathcal O}_{X\times \Aa^1}(1)=\text{sheaf
generated by }I\text{ in }{\mathcal O}_{X\times\Aa^1}(1).
$$
Choose a compactification $Y$ of $X\times\Aa^1$ to which
${\mathcal O}_{X\times\Aa^1}(1)$ extends to a line bundle
${\mathcal L}$ and let $\pi:B\to Y$ be a blow-up for
which $\pi^{-1}(J)={\mathcal O}_B(-E)$, where $E$ is the
exceptional divisor. Let $\pi^{-1}({\mathcal
L})={\mathcal O}_B(D)$, then we have
\begin{eqnarray}
\label{uno.nove}
e_{\bf r}(X)=D^d-(D-E)^d,
\end{eqnarray}
(\cite{Mo}, $\S 2$; \cite{Mu}, $\S 2$). One checks that
the property is independent of the choices involved. Let
$Z$ be the subscheme of $X\times\Aa^1$ defined by the
ideal sheaf $J$. According \cite{Fu} $\S 4$, one may
write the degree of contact in terms of Segre classes:
$$
e_{\bf
r}(X)=\sum_{k=1}^d(-1)^k\binom{d}{k}L^{d-k}s_k(Z,Y),
$$
(\cite{Fu}, Corollary $4.2.2$ and projection formula). Thus, when $Z$ is
set-theoretically a point $w_{\bf r}(m)$ is the
Hilbert-Samuel polynomial of $J$ as an ideal of
${\mathcal O}_{X\times\Aa^1}$ and $e_{\bf r}(X)$ the
multiplicity there (\cite{Mu}, $\S 2$ Examples).
\end{prg}

\begin{prg} \label{prg:due} The Chow point of $X$ is {\it semistable} with
respect to  integers $r_0\geq \cdots\geq r_N$ and a basis $l_0,\cdots,l_N$
of $E$ if and only if
\begin{eqnarray}
\label{uno.dieci}
E_\infty(X) \leq \frac{1}{n+1}\sum_{i=1}^Nr_i
\end{eqnarray}
(\cite{Mu}, Theorem $2.9$). If this properties is
satisfied for all bases and weights as above, the
Hilbert-Mumford theorem implies that the Chow point of
$X$ is semistable with respect to the action of $SL(E)$.
It is needless to say that the opposite of semistable is
{\it unstable}.
\end{prg}

\begin{prg} \label{prg:cinque} We finish this section with an important
remark. Schmidt Subspace Theorem (\cite{Sc}, Theorem
$V.1.D'$) corresponds to \thmref{th:uno} when $X$ is the
full projective space $\Pp^N$. We would like that
\thmref{th:uno} for an arbitrary projective variety
$X\subseteq \Pp^N$ is not a simple corollary of Schmidt
Subspace Theorem for $\Pp^N$. Paragraph
\prgref{prg:due} may help us with this task. It shows
that the possible varieties where this can be avoided
are those for which the Chow point of
$X\times\cdots\times X$ ($\# \Sigma$ times) is unstable.
\end{prg}

\section{Local Rings}

\begin{prg} Let $K\subseteq L$ be number fields, $\Sigma$ a finite
set of places of $L$ containing all infinite places. Let
$E$ be a vector space of rank $N+1$ over $K$ and $E^\vee$
$=$ Hom $(E,$ $K)$. Consider a closed subvariety
$X\subseteq
\Pp(E^\vee)$ in the projective space of lines of
$E^\vee$, of dimension of $d$. Let $P\in X(L)$ be
a closed point, and consider a basis $l_0,\cdots,l_N$ of
$E\otimes_KL$ so that $P=(1,0,\cdots ,0)$ with respect to
the coordinates in $\Pp(E^\vee)$ defined by this basis.
For $v\in\Sigma$ we define weights $r_{0,v}=k_v$,
$r_{1,v}=\cdots=r_{v,N}=0$.

\begin{prop}
Let $P\in X(L)$ be a closed point and for $v\in\Sigma$
let $d_v>0$ real such that
$$
\text{mult}_PX\sum_{v\in\Sigma}\frac{k_v}{d_v}>(\dim X+1)\deg X.
$$
Then the points ${\bf x}\in X(K)$ with
$$
\log\frac{|l_0({\bf x})|_v}{|{\bf x}|_v} \leq
-\frac{k_v}{s_v}h({\bf x}),\ v\in \Sigma,
$$
lye in finitely many subvarieties of $X$ of bounded
degree
\end{prop}

\proof For each $v\in\Sigma$ the ideal $J\cdot{\mathcal O}_{X\times
\Aa^1}(1)$ of
\propref{pr:tre} is generated by $\{t^{k_v}l_0,$ $l_1,$
$\cdots,$ $l_N\}$. Since $\{l_1,$ $\cdots,$ $l_N\}$
generate the maximal ideal ${\mathfrak m}_{P,X}$ and
$l_0$ is a unit at $P$, $J=(t^{k_v},{\mathfrak
m}_{P,X}){\mathcal O}_{X\times
\Aa^1}(1)$. Hence
$$
e_{\bf r}(X)=k_v\cdot\text{mult}_{(0,P)}(X\times
\Aa^1)=k_v\cdot\text{mult}_PX.
$$
The claim follows by the linearity of the expected value
in the weights, \formref{uno.cinque} and
\thmref{th:uno}.\qed
\end{prg}

\begin{prg} The attempt to generalize this theorem leads
to a numerical measure of the degree of singularity of a
point. Let $R$ be a local ring of dimension $r$ and $m$ a
positive integer, then the $m$-{\it th flat multiplicity} $e_m(R)$ of R
is defined by
\begin{eqnarray*}
e_0(R) & = & \sup\{\frac{e(I)}{r!col(I)}:\ I\text{ of
finite colength in }R\}\\ e_m(R) & = &
e_0(R[[t_0,\cdots,t_m]]),
\end{eqnarray*}
where $e(I)$ denotes the multiplicity of the ideal $I$.
Further, a local ring is called {\it semistable} if $e_1(R)=1$.
Let $X$ be a projective scheme. Let $L$ be an ample line
bundle on $X$. For a sufficiently large $n$ let
$$
\phi_n:X\to\Pp^{h^0(X,L^{\otimes n})-1},
$$
be the embedding defined by $L^{\otimes n}$. Suppose that
there exists a point $P$ on $X$ such that the local ring
${\mathcal O}_{P,X}$ is unstable. Then for every positive
iteger $n$, there exists a positive integer $m>n$ such
that the Chow point and the Hilbert point corresponding
to $\phi_n(X)$ are unstable under the natural action of
$SL(h^0(X,L^{\otimes m}))$ (\cite{Mu}, Propostion $3.12$,
\cite{Sh}, Proposition $1.3$). Little is known on the
semistability of local rings. See however \cite{Sh} and
\cite{Mu}, $\S 3$.
\end{prg}

\begin{prg} Let $K\subseteq L$ be number fields and $\Sigma$
a finite set of places of $L$ containing all infinte
places. Let $X$ be a projective scheme define over $K$
and $\Lb$ an ample line bundle on $X$. Suppose that there
is a point closed $P$ on $X(L)$ such that the local ring
$\Ob_{P,X}$ is unstable. By \cite{Mu}, Lemma $3.6$ there
exists a sequence of ideals of finite colenght
$$
I_0\subseteq I_1\subseteq\cdots\subseteq
I_N=\Ob_{P,X}=I_{N+1}=\cdots
$$
such that if $I$ is the ideal $\oplus_{i\geq
0}I_it^i\subseteq \Ob_{P,X}[[t]]$, then
\begin{eqnarray}
\label{localuno}
e(I)=(1+\varepsilon)(\dim X+1)!col(I),
\end{eqnarray}
where $\varepsilon>0$. By Riemann-Roch on $X$ we can
choose $n$ sufficiently large so that
\begin{enumerate}
\item $\Lb^{\otimes n}$ is very ample,
\item the map $S^mH^0(X,\Lb^{\otimes n})\to
H^0(X,\Lb^{\otimes nm})$ is surjective for all $m\geq 1$,
\item the map $\phi_m:H^0(X,\Lb^{\otimes mn})\to$ $\Ob_{P,X}/I_0^m$ is
surjective for all $m\geq 1$.
\end{enumerate}
The vector space $H^0(X,\Lb^{\otimes n})$ has the induced
filtration defined by
$$
F^j=\phi_i^{-1}(I_j/I_0),\ i=0,\cdots ,N.
$$
For each $v\in \Sigma$ choose a basis
$l_{v,0},\cdots,l_{v,N}$ of $H^0(X,\Lb^{\otimes n})$
compatible with the filtration. Assign weights as in
\formref{due:tre},
\begin{eqnarray}
\label{tre:uno}
r_{v,i}=\min\{ j:\ l_{v,i}\in F^j\}.
\end{eqnarray}
\begin{thm}
\label{th:cinque}
For $v\in\Sigma$ let $d_v>0$ be reals such that
\begin{eqnarray}
\label{localdue}
\sum_{v\in \Sigma}\sum_{i=0}^N\frac{r_{v,i}}{s_v}>
\frac{\deg X}{(\dim X+1)!}.
\end{eqnarray}
Then the points ${\bf x}\in X(K)$ with
\begin{eqnarray*}
\log(\frac{|l_{v,i}({\bf x})|_v}{|{\bf x}|_v}) \leq
-\frac{r_{v,i}}{d_v}h({\bf x}),\ v\in \Sigma,\ i=0,\cdots,N\\
\end{eqnarray*}
are contained in finitely many subvarieties of $X$ of
bounded degree.
\end{thm}

\proof Since the colength of $I$ corresponds to the sum of
weights(\cite{Sh}, p. $334$), the theorem follows by
\formref{localuno},\formref{localdue} and
\thmref{th:uno}.\qed
\begin{remark}
According to \cite{Sh} Appendix A, and \ref{prg:cinque},
\thmref{th:cinque} is stronger than the Schmidt
Subspace Theorem. For several examples of destabilizing
weighted flags of unstable two-dimensional local rings,
we refer to \cite{Sh}, $\S\S 4,5,6$.
\end{remark}
\end{prg}

\section{Ruled Surfaces}
\begin{prg}\label{prg:tre} The Steiner surface $X\subseteq \Pp^4$
is the closure of the image of the map
$\psi:\Pp^2\to\Pp^4$ defined by
$$
(x,y,z)\mapsto(xz,yz,x^2,xy,y^2)=(l_0,\cdots,l_4).
$$
This surface is the blown up at the point $P=(0,0,1)$ and
embedded by the system of conics pasing through $P$. Its
degree is three, the number of free points of
intersection of thwo such conics, and it is a rational
ruled surface of type $\Ff_1$, ruled by the pencil of
lines passing through $P$. Let us consider the flag given
by choosing the linear forms $l_i$ as above and with
weights $r_0=r_1=k$ for some $k>0$ and $r_2=r_3=r_4=0$.
Then then a computation shows that $e_{\bf r}(X)=4k$
(\cite{Mo}, Example $3.6$). \thmref{th:uno} implies:
\begin{prop}
Let $K\subseteq L$ be number fields, $\Sigma$ a finite
set of places of $L$ containg all infinite places and for
each $v\in \Sigma$ non-negative integers $k_v$ and
positive rationals $d_v$ with
$$
\sum_{v\in \Sigma}\frac{k_v}{d_v}>\frac{9}{4}.
$$
The solutions ${\bf x}\in X(K)$ of the inequalities
\begin{eqnarray*}
\log(\frac{|l_i({\bf x})|_v}{|{\bf x}|_v}) \leq
-\frac{k_v}{d_v}h({\bf x}),\ v\in \Sigma,\ i=0,1\\
\end{eqnarray*}
are contained in finitely many curves of $X$.
\end{prop}
\begin{remark}
If we
consider the same problem for ${\bf x}\in
\Pp^4(K)$, then Schmidt Subspace Theorem implies
that its solutions are contained in finitely many
subspaces of $\Pp^4$, hence in finitely many curves of
$X$, if the condition $\sum_{v\in S}(k_v/d_v)>5/2(>9/4)$
holds true.
\end{remark}
\end{prg}

\begin{prg} A good reference for the general properties of ruled
surfaces is \cite{Ha}, $V.2$. Fix a smooth curve $C$ of
genus $g$ and a geometrically ruled surface $p:R\to C$.
For one section $\sigma:C\to R$ of $p$ we refer to the
divisor $S$ on $R$ which is the image of $\sigma$. Fix
one section, say $S$, and let $f$ denote the numerical
equivalence class of a fibre of $p$. When the parity of
$S^2$, which is independent of $S$, is even. A convenient
basis of $Num(R)$ is given by $f$ and the element $G$
determined by $G^2=0$ and $G\cdot f=1$. For suraces of
odd parity such a $G$ can be found in $Num(R)\otimes
\Qq$. We will by abuse of language consider $G$ as a divisors on
$R$. This is no restriction since we will be utilizing
only the numerical properties of $G$.
\end{prg}

\begin{prg} Fix a ruled surface $R$, a very ample divisor $D\sim aG+bf$
on $R$, and a section $S$ such that $D-aS$ is effective.
Let $\Pp^N=\Pp(H^0(R,D)^{\vee})$ and let $F^j$ $=\{$ $h$
$\in$ $H^0(R,D),$  $h\in H^0(R, D-jS)\}$. Since a section
in $D$ can vanish to order at most $a$ on $S$,
\begin{eqnarray}
\label{due:due}
H^0(R,D)=F^0\supset F^1\supset\cdots\supset F^a\supset
F^{a+1}=\{0\}.
\end{eqnarray}
Choose a basis $l_0,\cdots ,l_N$ of $H^0(R,D)$ compatible
with this flag, and fix weights $r_i$, $i=0,\cdots ,N$ by
the condition that
\begin{eqnarray}
\label{due:tre}
r_i=a-j\Longleftrightarrow x_i\in F^j\setminus F^{j+1}.
\end{eqnarray}
That is, $r_i$ equals $a$ minus the order to which $l_i$
vanishes along $S$, for $i=0,\cdots ,N$. Clearly, the
$r_i$'s decrease to zero. This construction generalize
\prgref{prg:tre} of the Steiner surface in $\Pp^4$. In that
construction we have $C=\Pp^1$, $E={\mathcal
O}_{\Pp^1}(1)\otimes {\mathcal O}_{\Pp^1}(2)$ and
$R=\Pp(E)$. Further, $D={\mathcal O}_{\Pp(E)}(1)\sim
G+3/2f$ and $S$, the exceptional divisor, which is
numerical equivalent to $G-1/2f$ and which is the section
associated to the bundle ${\mathcal O}_{\Pp^1}(2)$ in
$E$, we get the flag of that example. By \cite{Mo} Theorem $6.5$ we
obtain a non-trivial generalisation of Schmidt Subspace Theorem for
this embedding:
\begin{thm}
\label{th:tre}
Let $K\subseteq L$ be number fields and $\Sigma$ a finite
set of places of $L$ containing all infinite places. Let
$S$ be a section of $R$. For each $v\in\Sigma$ let
$l_{v,0},\cdots ,l_{v,N}$ a basis compatible with the
filtration
\formref{due:due} associated to $S$ and let
$r_{v,0},\cdots,r_{v,N}$ be the corresponding weights
\formref{due:tre}. Suppose that $H^i(R,D-jS)=0$ for $i>0$
and $0\leq j\leq a$. Further, for each $v\in \Sigma$ let
$d_v>0$ reals such that
\begin{eqnarray}
\label{due:uno}
(3a^2D\cdot S-a^3S^2)\sum_{v\in\Sigma}\frac{1}{d_v}>3D^2.
\end{eqnarray}
Then the points ${\bf x}\in R(K)$ with
$$
\log\frac{|l_i({\bf x})|_v}{|{\bf x}|_v} \leq
-\frac{r_{v,i}}{d_v}h({\bf x}),\ v\in \Sigma,\ i=0,\cdots,N
$$
lie in finitely many curves of $R$ of bounded degree.
\end{thm}
\proof This follows by \thmref{th:uno} and \cite{Mo}
Proposition $6.2$.\qed
\end{prg}

\begin{prg} If $p:E\to C$ is a rank $2$ vector bundle on $C$ such
that $\Pp(E)\cong R$ we say that $E$ represent $R$. Such
$E$ always exists and $R$ determines $E$ up to tensoring
with a line bundle $L$ on $C$. Let $L$ be a line
subbundle of $E$ of maximal degree. The number $\deg
E-2\deg L$ is independent of the choice of $E$. Then we
say that the ruled surface is {\it bundle semistable} if
$\deg E\geq 2\deg L$ and {\it bundle unstable} when $\deg
E<2\deg L$. In the latter case there exists a unique
section $S$ of negative selfintersection.
\begin{thm}
\label{th:quattro}
Let $K\subseteq L$ be number fields, $\Sigma$ a finite
set of places of $L$ containing all infinite places.
Suppose $R$ is bundle unstable. Let $D\sim aG+bf$ be a
very ample on $R$. Let $S$ be the unique section of $R$
of negative selfintersection. For $v\in\Sigma$ let
$l_{v,0}$, $\cdots,$ $l_{v,N}$ a basis compatible with
the filtration
\formref{due:due} associated to $S$ and let
$r_{v,0},\cdots,r_{v,N}$ be the corresponding weights
\formref{due:tre}. Further, for each $v\in \Sigma$
let $d_v>0$ reals such that
\begin{eqnarray}
\label{fo:dueuno} b +\frac{a}{2}(\deg
E-2\deg L) >  2g-2,
\end{eqnarray}
and
\begin{eqnarray}
\label{fo:duedue} \sum_{v\in\Sigma}\frac{1}{d_v}>
\frac{3b}{a(b+g-1)}.
\end{eqnarray}
Then the points with ${\bf x}\in R(K)$ with
$$
\log\frac{|l_i({\bf x})|_v}{|{\bf x}|_v} \leq
-\frac{r_{v,i}}{d_v}h({\bf x}),\ v\in \Sigma,\ i=0,\cdots,N
$$
lie in finitely many curves of $R$ of bounded degree.
\end{thm}

\proof Since $\deg
E<2\deg L$ and $b$ is bounded below by
\formref{fo:dueuno}, $H^i(R,D-jS)=0$ for $i>0$ and
$0\leq j\leq a$. (\cite{Mo}, Remark after Proposition
$6.2$). The degree of $R$ is $D^2=2ab$. According to
\cite{Mo} $(6.7)$ we have $3a^2D\cdot
S-a^3S^2=\frac{a^3}{2}l_R+3a^2b$. The inequalities
\formref{fo:dueuno} and \formref{fo:duedue} imply
\formref{due:uno} and \thmref{th:tre} concludes the theorem.
\qed
\end{prg}

\begin{prg} Little is known of the unstability of higher dimensionals varieties. Each new information on this problem has always interesting consequences. For instance, Vojta's Conjectures on the ruled threefold defined by the equation
$a^4x+b^4y+c^4z=0$ in $\Pp^2\times\Pp^2$ imply $abc$-conjecture.
\end{prg}

\section{Blow-up}

\begin{prg} Let $K\subseteq L$ be number fields, $\Sigma$ a finite
set of places of $L$ containing all infinite places. Let
$X$ be a normal variety defined over $K$ of dimension
$n$, $H$ a Cartier divisor on $X$, $E$ an effetcive
Cartier divisor and fix an integer $s>0$. Assume that for
some sufficiently large $r$, larger than $r_0\geq s$ say,
$rH-sE$ is very ample. Put
$W_{i,j}=H^0(X\otimes_KL,iH-jE)$ for $i,j\in
\Zz$. Then
$W_{r,s}$ has a filtration $W_{r,s}\supset
W_{r,s+1}\supset
\cdots\supset W_{r,r}$. For $v\in\Sigma$ let $l_{v,0}$,
$l_{v,1}$, $\cdots$, $l_{v,N}$ be a basis of $W_{r,s}$
compatible with the filtration. For these bases we
consider weights $r_{v,i}=r-j$ if $l_{v,i}\in W_{r,j}$
and $l_{v,i}\not\in W_{r,j+1}$, $v\in\Sigma$.

\begin{thm}
\label{th:sette} Let us suppose that for any $i>0$
$$
h^i(X\otimes_KL,rH-jE)=O((r+j)^{n-2})\text{ on
}\{(r,j)\in\Zz^2|\ s\leq j\leq r\}.
$$
For $v\in\Sigma$ let $d_v$ be positive reals with
\begin{eqnarray}
\sum_{i=0}^nH^{n-i}(-E)^i\binom{n+1}{i+1}((i+1)(r-s)r^{n-i}
s^i-r^{n+1}+r^{n-i}s^{i+1})\sum_{v\in\Sigma}\frac{1}{d_v}>
(n+1)(rH-sE)^n.\\ \label{Hn-i-Ei}
\end{eqnarray}
Then the points ${\bf x}\in X(K)$ with
\begin{eqnarray*}
\log(\frac{|l_{v,i}({\bf x})|_v}{|{\bf x}|_v}) \leq
-\frac{r_{v,i}}{d_v}h({\bf x}),\ v\in \Sigma,\ i=0,\cdots,N\\
\end{eqnarray*}
are contained in finitely many subvarieties of $X$ of
bounded degree.
\end{thm}
\proof For each $v\in \Sigma$ the filtrations constructed above correspond to the
filtration of \cite{Is}, Lemma $3$. In {\it loc. cit.} it
is proven that,for each $v\in\Sigma$ the left hand side
of \formref{Hn-i-Ei} is a lower bound for the degree of
contact. The assertion we want follows then by
\thmref{th:uno}.\qed

\begin{remark} Lemma $3$ of  \cite{Is} confirms that this theorem
is stronger than Schmidt Subspace Theorem restricted to
$X$. Further, one should observe that all intersection
products of \formref{Hn-i-Ei} can be written down by
Segre classes.
\end{remark}
\end{prg}

\begin{prg} In this section we follow \cite{Is}, $\S 2$.
Let $V$ be a normal projective variety of dimension
$n-1\geq 1$ over $K$, $D$ an ample divisor on $V$ and
$Y_t\hookrightarrow\Pp^{N(t)+1}$ the projective cone over
$\Phi_{|tD|}:V\hookrightarrow\Pp^{N(t)}$. Let $f_t:X_t\to
Y_t$ be the blow-up with center the vertex. Put
${\mathcal O}_{Y_t}(H_{0,t})={\mathcal
O}_{\Pp^{N(t)+1}}(1)|_{Y_t}$, $H_t=f_t^*H_{0,t}$ and
$E_t$ the exceptional divisor of $f_t$.
\begin{prop}
\label{pr:blowup} Take a large integer $t$ in such
a way that
\begin{eqnarray}
\label{hivmtl}
H^i(V\otimes_KL,mt{\mathcal L})=0,\text{ for any
}m>0\text{ and }i>0.
\end{eqnarray}
Further, assume that
\begin{eqnarray}
\label{tndeglv}
(n+1)s^n(s-r)+r^{n+1}-s^{n+1}>\frac{n+1}{\#\Sigma}(r^n-s^n).
\end{eqnarray}
Then the result of \thmref{th:sette} holds.
\end{prop}
\proof Condition \formref{hivmtl} implies the vanishing of
higher cohomologies for each $rH_t-jE_t$, $j=1,\cdots,r$
(\cite{Is}, Proposition $4$). Moreover, all intersection
numbers $H_t^{n-i}(-E_t)^i$ are zero except for $i=0,$
$n$ where $H_t^{n}=-(-E_t)^{n}$. Thus,
\formref{tndeglv} implies \formref{Hn-i-Ei} and the
proposition follows.\qed
\begin{remark}
The projective cone over a non-singular conic in $\Pp^2$
and a non-singular quadratic surface in $\Pp^3$ satisfy
the condition of \propref{pr:blowup} for each $t\geq 1$
(\cite{Is}, Example $6$.)
\end{remark}
\end{prg}
\begin{prg} We consider here resolutions of singularities. Let
$f:X\to Y$ be the blow-up of a normal projective variety
$Y$ whose center is a subscheme $Z$ not necessarily
irreducible or reduced. Assume $X$ is non-singular of
dimension $n$. Let $E$ be the exceptional divisor. Put
$H=f^*H_0$ for an ample divisor $H_0$ on $Y$, and fix an
integer $s\geq 0$.
\begin{prop}
\label{pr:ressing}
Take a sufficiently large integer $t$ such that $tH-E$ is
ample on $X$ and $tH_0$ is very ample on $Y$. Assume
\formref{Hn-i-Ei}, then the assertion of
\thmref{th:sette} holds.
\end{prop}
\proof The vanishing for the higher
cohomology spaces is assured by the proof of \cite{Is},
Proposition $8$.\qed
\end{prg}

\section{Elliptic surfaces}

\begin{prg} Following \cite{Mi}, we call
a flat proper map of $\overline{\Qq}$-schemes
$p:X\to\Pp^1$ a rational Weierstrass fibration if $X$ is
reduced and irreducible rational surface over
$\overline{\Qq}$, every geometric fibre of $p$ is an
irreducible curve of genus $1$, and a section $s:\Pp^1\to
X$ is given, not passing through the nodes or cusps of
the fibres. Moreover, we will assume that $X$ is normal,
and that the generic fibre of $p$ is smooth. In this
case, we may resolve the singularities of $X$ and obtain
an elliptic surface $\overline{p}:\overline{X}\to\Pp^1$
(with section) which we call the induced elliptic
surface. In this situation one may represent $X$ in
Weierstrass form by the equation
$$
y^2=x^3+A(t)x+B(t),
$$
where $A$ is a quartic and $B$ is a sextic polynomial in
the parameter $t$ of $\Pp^1$.
\end{prg}

\begin{prg} Let us fix the following situation: $\overline{X}$
is a minimal rational elliptic surface with section and a
fibre of type $IV^*$ (in Kodaira classification), $X$ the
associated Weierstrass fibration, $S$ the given section
of $\overline{X}$, $f$ the numerical class of a fibre of
$\overline{X}$, $E$ the unique rational component of the
singular $IV^*$ fibre having multiplicity $3$,
$D=3kS+6kf$ a very ample divisor on $X$ for some positive
integer $k$. We have a filtration on $H^0(\overline{X},
{\mathcal O}(D))$ defined by
\begin{eqnarray}
H^0(\overline{X}, {\mathcal O}(D))\supset
H^0(\overline{X}, {\mathcal O}(D-E))\supset\cdots\supset
H^0(\overline{X}, {\mathcal O}(D-18kE))\supset
\{0\}.\nonumber\\
\label{elluno}
\end{eqnarray}
Let $K$ be a number field defining $\overline{X}$ as a
geometrically irreducible variety, and $L$ a finite
extension. Let $\Sigma$ a finite set of places of $L$
containing all infinite places. For each $v\in\Sigma$, we
choose a basis $l_{v,0},\cdots,l_{v,n}$ compatible with
this filtration, and define weights $r_{v,i}=18k-j$ if
$l_{v,i}$ is in $H^0(\overline{X}, {\mathcal O}(D-j E))$
but not in $H^0(\overline{X}, {\mathcal O}(D-(j+1)E))$,
$i=0,\cdots,n$.
\begin{thm}
\label{th:otto}
For $v\in\Sigma$ let $d_v$ be positive reals such that
$$
\sum_{v\in\Sigma}\frac{k}{d_v}>\frac{3}{25}.
$$
Then all ${\bf x}\in X(K)$ with
\begin{eqnarray*}
\log(\frac{|l_{v,i}({\bf x})|_v}{|{\bf x}|_v}) \leq
-\frac{r_{v,i}}{d_v}h({\bf x}),\ v\in \Sigma,\ i=0,\cdots,n\\
\end{eqnarray*}
are contained in finitely many curves of $X$ of bounded
degree.
\end{thm}
\proof According to \cite{Mi}, $(6.7)$ we know that
the degree of contact of the filtration \formref{elluno}
is bounded below by $675k^3$. On the other hand, the
degree of $X$ with respect to $D$ equals $27k^2$
(\cite{Mi}, $(4.1)$). Hence,
$$
\frac{e_{\bf r}(X)}{(\dim X+1)\deg X}\geq
\frac{675k^3}{3\cdot 27k^2}=\frac{25}{3}k.
$$
\thmref{th:uno} confirms the the theorem.\qed

\begin{remark} Since $D$ is very ample, by Riemann-Roch we can
compute the dimension $n+1$ of the space
$H^0(\overline{X}\otimes_KL, {\mathcal O}(D))$ which
equals $(27k^2+3k+2)/2$ (\cite{Mi}, $(4.2)$). For each
$v\in\Sigma$
$$
\sum_{i=0}^nr_{v,i}=\frac{225}{2}k^3+9k^2+\frac{23}{2}k,
$$
(\cite{Mi}, $(5.2)$). Moreover, for all $k\geq 1$ we have
$$
\frac{25}{3}k>\frac{1}{\frac{27k^2+3k+2}{2}}(\frac{225}{2}k^3+
9k^2+\frac{23}{2}k).
$$
This means that Schmidt Subspace Theorem restricted to
$X$ is weaker than \thmref{th:otto}
\end{remark}

\begin{remark} The filtration \formref{elluno} is clearly given by the
order of vanishing for sections of $D$ along the curve
$E$. Calculations of the same type can be carried through
for Weierstrass fibrations with singularities of type
$II^*$ and $III^*$. The filtrations used in these latter
cases are also given by the order of vanishing along the
unique curve of maximum multiplicity in the singular
fibre in question. This curve is the multiplicity six
component in the $III^*$ fibre and the multiplicity four
component in the $II^*$ fibre.
\end{remark}
\end{prg}

\end{document}